\documentclass[11pt]{article}
\usepackage{amsmath,amsfonts,amsthm,epsf}
\usepackage{graphicx,subfigure}
\usepackage[export]{adjustbox}
\usepackage{epstopdf}

\usepackage{hyperref}
\usepackage{enumerate}
\usepackage{verbatim}
\usepackage{tcolorbox}

\usepackage{multirow}

\newcommand{\ve}{\varepsilon}

\RequirePackage{natbib}

\usepackage[margin=2cm]{geometry}

\newcommand{\ignore}[1]{}

\newtheorem{theorem}{Theorem}
\newtheorem{proposition}{Proposition}
\newtheorem{lemma}{Lemma}

\begin{document}

\title{Optimal minimax random designs for weighted least squares estimators}
\author{
David Azriel, Technion – Israel Institute of Technology
}
\date{}
\maketitle
\begin{abstract}	
	This work studies an experimental design problem where {the values of a predictor variable, denoted by $x$}, are to be determined with the goal of estimating a function $m(x)$, which is observed with noise. A linear model is fitted to $m(x)$ but it is not assumed that the model is correctly specified. It follows that the quantity of interest is the best linear approximation of $m(x)$, which is denoted by $\ell(x)$. It is shown that in this framework the ordinary least squares estimator typically leads to an inconsistent estimation of $\ell(x)$, and rather weighted least squares should be considered. 
	An asymptotic minimax criterion is formulated for this estimator, and a design that minimizes the criterion is constructed. An important feature of this problem is that the $x$'s should be random, rather than fixed. Otherwise, the minimax risk is infinite. It is shown that the optimal random minimax design is different from its deterministic counterpart, which was studied previously, and a simulation study indicates that it generally performs better when $m(x)$ is a quadratic or a cubic function. Another finding is that when the variance of the noise goes to infinity, the random and deterministic minimax designs coincide. 
	The results are illustrated for polynomial regression models and the general case is also discussed.    
\end{abstract}

\noindent {\bf Key words:}	Experimental design, I-optimality, Robust regression.

\section{Introduction}

A fundamental problem in experimental design is to {determine the values of a predictor variable $x$ } in order to estimate a function $m(x)$, which is observed with noise. Specifically, consider a pair $(X,Y)$ of univariate random variables and suppose that the quantity of interest is $m(x)=E(Y|X=x)$. Assume also that $X$ can assume values in the interval $[-1,1]$.  The experimenter selects $X_1,\ldots,X_n$, which results in the corresponding observations $Y_1,\ldots,Y_n$. Then, based on $(X_1,Y_1),\ldots,(X_n,Y_n)$ an estimator $\hat{m}(x)$ is considered, and the question of interest is how to select $X_1,\ldots,X_n$ in order to minimize  $\int_{-1}^1 E\{\hat{m}(x)-m(x)\}^2 dx$. This problem is referred to as Q-optimality in \citet{Fedorov}; other authors  \citep[e.g.,][]{Wiens2015} use the term I-optimality.

When a polynomial of degree $K$ is used to model $m$, i.e., $m(x)=\sum_{k=0}^K b_k x^{k}$ for some unknown coefficients $b_0,\ldots,b_k$, which are estimated using least squares, a classical result shows that the support of the optimal sampling measure contains $K+1$ points. These points can be computed as well as the optimal weights \citep{Kiefer_Wolfowitz}.  For the special case where $m$ is a linear function, i.e., $m(x)=b_0 + b_1 x$, the optimal design is to sample $X=\pm1$ with equal weights.

The above designs do not allow one to check the model assumptions and therefore are not robust to model misspecification. In order to address this problem, \citet{box1959basis} considered the situation where a polynomial model of degree $K$ is assumed, but the true function $m$ is a polynomial of higher degree. This seminal paper started a line of work where model misspecification is accounted for in order to compute optimal designs. A review of this literature can be found in \citet{Wiens2015}; see also the recent work of \citet{Waite}. A major theme in this literature is that the bias term is more important than the variance, and therefore one can focus on designs that minimize the bias. This idea can be summarized in the following quote from Box and Draper:
\begin{quote}
	``We are led therefore to a somewhat remarkable conclusion. This is that the optimal design in a ``typical" situation in which the influence of bias and variance are equal is very nearly the same as that obtained when variance is ignored completely and the experiment is designed to minimize bias alone.'' (page 630)
\end{quote}
However, as I argue below, this is a result of using the ordinary least squares estimator, which is inconsistent in this setting to the parameter of interest -- the best linear approximation. When using a consistent estimator, namely the weighted least squares estimator, the bias term is negligible compared to the variance.  

{Another important consideration is the randomness of the design - an idea that goes back at least to \citet{Fisher1926} and has a long and rich history in  experimental design. It still plays a dominant role in current online A/B testing \citep{Kohavi2017} as well as in recent clinical trials \citep{Knoll2021}.
	\citet{Cox2009} reviews this large body of literature and mentions three objectives of randomization: the avoidance of bias, error estimation, and the provision of exact tests. In the context of the current work, the first of these objectives is essential. Indeed, \citet{Wiens1992} shows that deterministic designs have infinite minimax risk as the resulting bias can be arbitrarily large. 
	\citet{Waite} study a random experimental design strategy in the above setting and construct a specific random design that has a smaller minimax error than its deterministic counterpart.
	In the current work, the random experimental design approach is adopted as detailed below.}    

To fix ideas, a simple linear approximation is now considered. 
The best linear approximation of $m(x)$ with respect to $X\in U[-1,1]$ is denoted by $\ell(x)=\beta_0 + \beta_1 x$, where
\begin{equation}\label{eq:m_ell}
	\left( \begin{array}{c}{\beta_{0}}\\{\beta_{1}}\end{array}\right)= {Q}^{-1} \left( \begin{array}{c} \int_{-1}^1 \frac{1}{2} m(x) dx \\ { \int_{-1}^1 \frac{x}{2} m(x) dx } \end{array}\right) \text{ and }{Q}=  \left( \begin{array}{cc} 1 &\int_{-1}^1 \frac{x}{2} dx\\ \int_{-1}^1 \frac{x}{2} dx & \int_{-1}^1 \frac{x^2}{2} dx \end{array}\right)=
	\left( \begin{array}{cc} 1 &0\\ 0& 1/3 \end{array}\right);
\end{equation}
see Section 2.25 of \citet{Hansen}. For a given estimator of $\ell(x)$, which is denoted by $\bar{\ell}(x)$ and assumed to be a linear function of $x$, the experimenter seeks to minimize $\int_{-1}^1 E\{\bar{\ell}(x)-m(x)\}^2 dx$. {The definition of $\ell(x)$ in \eqref{eq:m_ell} implies the orthogonality property of $\ell(x)$,  i.e., $\displaystyle{\int_{-1}^1 \binom{1}{x}}\{\ell(x)-m(x)\}dx=\displaystyle{\binom{0}{0}}$; see also Section 2.18 of \citet{Hansen}}. It follows that
\[
\int_{-1}^1E\{\bar{\ell}(x)-m(x)\}^2 dx =  \int_{-1}^1 E\{\bar{\ell}(x)-\ell(x)\}^2 dx+ \int_{-1}^1 \{{\ell}(x)-m(x)\}^2 dx.
\] 
The second term does not depend on the design and therefore can be ignored. Continuing with the first term,
\begin{equation}\label{eq:bar_ell}
	\int_{-1}^1 E\{\bar{\ell}(x)-\ell(x)\}^2 dx= \int_{-1}^1 {\rm var}\{\bar{\ell}(x)\} dx+ \int_{-1}^1 [E\{\bar{\ell}(x)\}-\ell(x)]^2 dx. 
\end{equation}
The first term is referred to as the variance and the second as bias in \cite{box1959basis}. As mentioned previously, Box and Draper found that the bias term dominates the variance. However, this depends on the estimator used, as explained below.  

{It is assumed from now on that the experimenter uses a random design strategy, that is $X_1,\ldots,X_n$ are sampled from a distribution on $[-1,1]$ with density $\pi(x)$ with respect to the Lebesgue measure. Thus, the design, $X_1,\ldots,X_n$, is random and the experimenter wishes to find the optimal $\pi$ according to a minimax criterion, which is defined below. }
It is assumed that $\pi(x)>0$ for all $x \in [-1,1]$; otherwise, the risk in \eqref{eq:bar_ell} could be infinite for various settings; see \citet{Wiens1992} and \citet{Heo2001}. The corresponding observations are $Y_1,\ldots,Y_n$. The ordinary least squares estimators, which are denoted by $\hat{\beta}_0,\hat{\beta}_1$, are a function of $(X_1,Y_1),\ldots,(X_n,Y_n)$ and the estimate of the best linear approximation is $\hat{\ell}(x)=\hat{\beta}_0 + \hat{\beta}_1 x$. 

\begin{figure}[h]
	\subfigure[$m(x)$, $\ell(x)$ and $\ell_\pi(x)$]{%
		\includegraphics[width=0.45\textwidth]{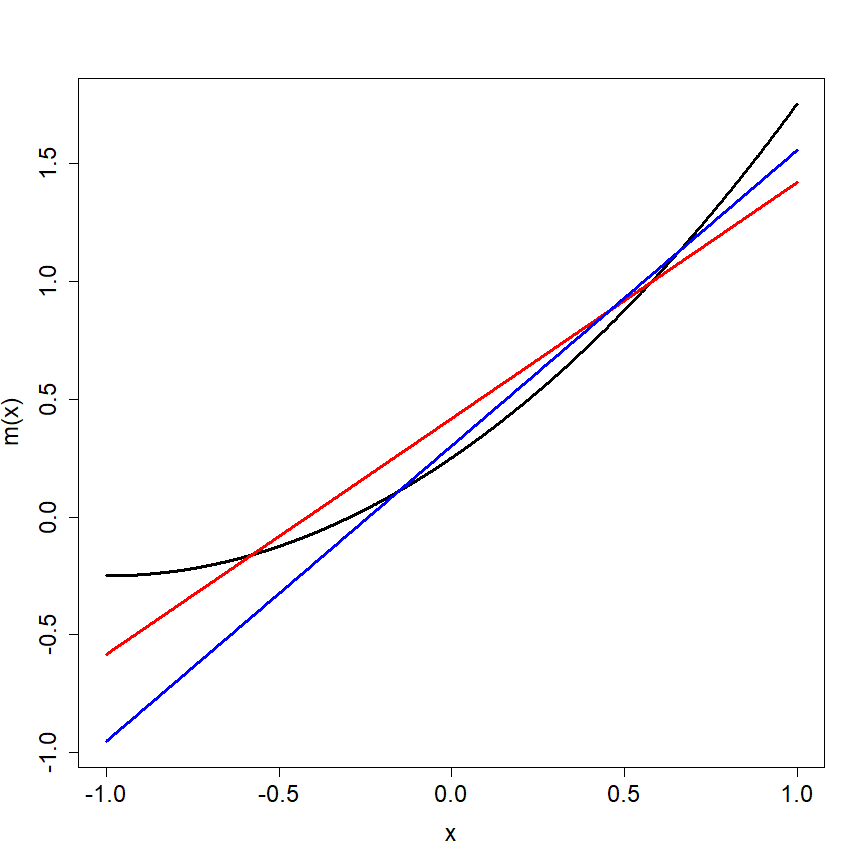}
	}
	\subfigure [Marginal distribtion of $x$]{%
		\includegraphics[width=0.45\textwidth]{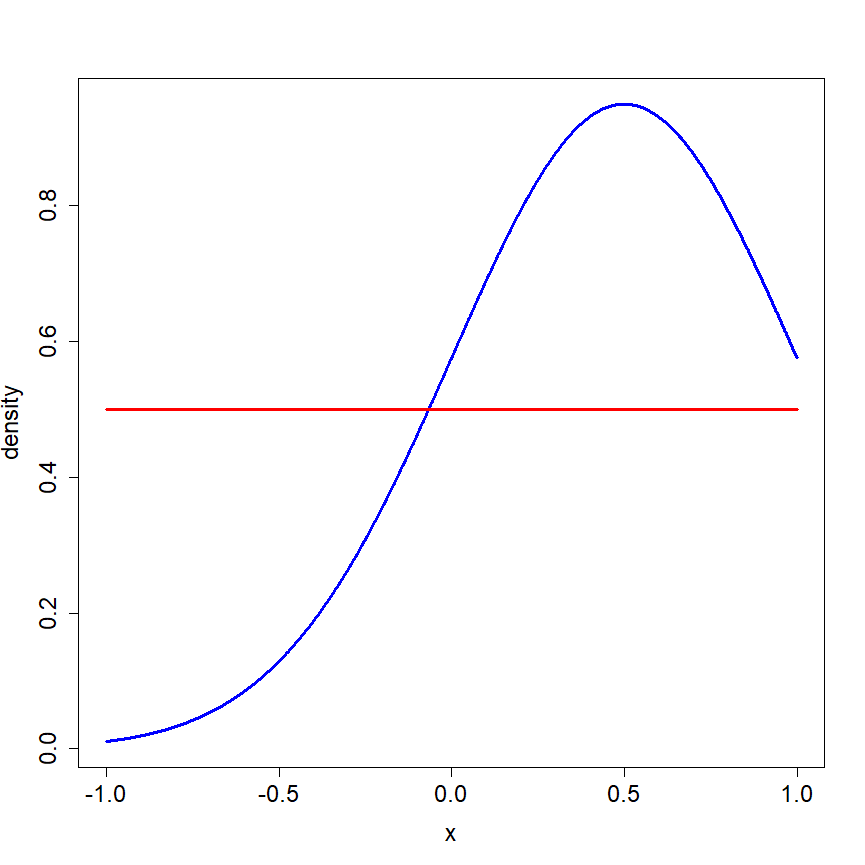}
	}
	\caption{\footnotesize Figure (a): the conditional expectation $m(x)$ (black) and the best linear approximation for two possible marginal distributions of $X$: uniform (red) and $\pi=$ truncated $N(1/2,1/4)$ (blue). Figure (b): The densities of these two distributions.}\label{fig:BLA}
\end{figure}

The important thing to notice here is that $\hat{\ell}(x)$ converges to the best linear approximation of $m(x)$ with respect to $\pi$, rather than $\ell(x)$, which is the best linear approximation when $X\in [-1,1]$. Specifically, the limit of $\hat{\ell}(x)$ is $\ell_\pi(x)=\beta_{0,\pi}+\beta_{1,\pi}x$, where
\[
\left( \begin{array}{c}{\beta_{0,\pi}}\\{\beta_{1,\pi}}\end{array}\right)= {Q}_\pi^{-1} \left( \begin{array}{c} \int_{-1}^1 \pi(x) m(x )dx \\ { \int_{-1}^1 x\pi(x) m(x) dx } \end{array}\right) \text{ and }{Q}_\pi=  \left( \begin{array}{cc} 1 &\int_{-1}^1 x \pi(x) dx\\ \int_{-1}^1 x \pi(x) dx & \int_{-1}^1 x^2 \pi(x) dx \end{array}\right).
\] 
{A subscript $\pi$ is used to denote quantities that depend on $\pi$. }
To illustrate this point, consider a specific setting where $m(x)=1/4+x/2 + x^2/4$ and $\pi$ is the density of truncated $N(1/2,1/4)$ on the interval $[-1,1]$. This setting is illustrated in Figure \ref{fig:BLA}. For the truncated normal, more weight is given for large $x$ and therefore the intercept of $\ell_\pi(x)$ is smaller and the slope is larger.

The above discussion implies that $\hat{\ell}(x)$ is an inconsistent estimator of $\ell(x)$ unless $\pi$ is the uniform density. This has two consequences: first, the bias term dominates the variance in \eqref{eq:bar_ell}, and second, asymptotically, the optimal design is the uniform design. These two findings were described in \citet{box1959basis} and in other papers; see Section 20.3 of \citet{Wiens2015} for a literature review. However, they are a result of using an inconsistent estimator, rather than of robustness considerations.  Below, a consistent estimator to $\ell(x)$ is formulated, and it is found, in contrast, that the bias term vanishes asymptotically and the uniform design is no longer optimal.   

In this framework, there is another reason to consider $X$ as random rather than fixed, besides the one that was mentioned above.
\citet{buja2019models} show that $X_1,\ldots,X_n$ is ancillary to $\beta_0,\beta_1$, if and only if $m(x)=\ell(x)$. This means that when the linear model $m(x)=\ell(x)$ is not assumed, the marginal distribution of $X$ matters, and the conditionality principle \citep[][page 38]{Cox_1974} does not apply, i.e., the $X$'s should not be conditioned upon. 

{The goal of this work is to find an optimal density $\pi$ that is optimal in a minimax sense. \citet{Nie2018} study the analogous problem in the active learning setting where the purpose is to select predictors from a given unlabeled dataset according to the minimax criterion. In particular, the problem addressed in Section 3.1 of  \citet{Nie2018} is similar to the one that is considered in the current work and the theoretical results are similar in both papers. The current work is also close to  \citet{Wiens1998} and \citet{Wiens2000}. A detailed comparison to these works is given below in Section \ref{sec:main}. 
}

\section{Problem formulation} \label{sec:prel}

A linear approximation to $m(x)=E(Y|X=x)$ is considered, using the vector of monomials of $x$, $\vec{x}=(1,x,\dots,x^K)$; general functions of $x$ are considered in Section 6 of the supplementary material. The coefficient of the best linear predictor is ${\beta}={Q}^{-1} \int_{-1}^{1}\frac{1}{2} \vec{x}^T m(x)dx$, where ${Q}= \int_{-1}^{1} \frac{1}{2} \vec{x}^T \vec{x} dx$. In this notation, ${\beta}$ and  $\vec{x}$ are column and row vectors, respectively. The best linear approximation is the function $\ell(x)=\vec{x} {\beta}$.

The experimenter chooses a density $\pi$ from which $X_1,\ldots,X_n$ are sampled. Then $Y_1,\ldots,Y_n$ are observed, where
\begin{equation}\label{eq:model}
	Y_i= m(X_i)+\ve_i,
\end{equation}
and the residuals $\ve_1,\ldots,\ve_n$ are independent and identically distributed satisfying $E(\ve_i^2|X_i)=\sigma^2$, i.e., homoscedasticity is assumed; this assumption is relaxed in Section \ref{sec:var} below.  
Let $\vec{{X}}_i=(1,X_i,\ldots,X_i^K)$ and define the regression residuals $e_i = Y_i - \vec{{X}}_i{\beta}$, $i=1,\ldots,n$. The vector ${\beta}$ is the best linear predictor with respect to $X \sim U[-1,1]$ and not to $X \sim \pi$. Hence, $E_\pi( \vec{{X}}_1 e_1)$ is not necessarily zero. However, when $\vec{{X}}_1$ is weighted by $1/\pi(X_i)$, then the expectation is zero, i.e., $E_\pi( \vec{{X}}_1 e_1/\pi(X_1))={0}$.

Let ${\mathbb{X}} \in \mathbb{R}^{n \times (K+1)}$ be the matrix with $\vec{{X}}_i$'s in the rows, ${\mathbb{Y}}=(Y_1,\ldots,Y_n)^T$ and let ${W}_\pi \in \mathbb{R}^{n \times n}$ be the diagonal matrix whose $i$-th diagonal element is $\frac{1/2}{\pi(X_i)}$. Define
\begin{equation} \label{eq:tilde_Q}
	\tilde{Q}_\pi = I({\cal E}_\pi) n {Q} +  I({\cal E}_\pi^C){\mathbb{X}}^T {W}_\pi {\mathbb{X}},
\end{equation}
where $I$ denotes the indicator function and
\begin{equation}\label{eq:calE}
	{\cal E}_\pi= \{ \lambda_{min}({\mathbb{X}}^T {W}_\pi {\mathbb{X}}/n) < \lambda_{min}({Q})/2  \},
\end{equation}
and where $\lambda_{min}$ denotes the smallest eigenvalue. A discussion about the estimator $\tilde{Q}_\pi$ is given below. Define also $\tilde{V}_\pi= {\mathbb{X}}^T {W}_\pi {\mathbb{Y}}$. 
The weighted least squares estimate is $\tilde{\beta}_\pi =\tilde{Q}_\pi^{-1} \tilde{V}_\pi$ and the estimate of the best linear approximation is $\tilde{\ell}_\pi(x)=\vec{x}\tilde{\beta}_\pi$.

The specific choice of $\tilde{Q}_\pi$ in \eqref{eq:tilde_Q} is now discussed. Under the setting considered here ${Q}$ is known, and hence one can use the estimator ${Q}^{-1} \tilde{V}_\pi/n$. However, in this context this estimator is typically less efficient; see \citet{Tarpey} and \citet{Cook}. Another option is the estimator $\tilde{\beta}_{\pi,0}= \left({\mathbb{X}}^T {W}_\pi {\mathbb{X}}\right)^{-1} {\mathbb{X}}^T {W}_\pi {\mathbb{Y}}$. This is an efficient estimator for ${\beta}$, but for finite $n$, $\lambda_{min}\left({\mathbb{X}}^T {W}_\pi {\mathbb{X}}/n\right)$ can be very small with positive probability and therefore the expectation of $\tilde{\beta}_{\pi,0}$ may not exist; see Section 4.7 of \citet{Hansen}. On the other hand,   
the event ${\cal E}_\pi$ occurs with exponentially small probability \citep[][Theorem 5.1]{Tropp} and therefore, in practice the difference between $\tilde{\beta}_\pi$ and $\tilde{\beta}_{\pi,0}$ is small. In particular, they have the same asymptotic distribution. The estimator $\tilde{\beta}_\pi$ is advantageous because $(\tilde{ Q}/n)^{-1}$ is bounded and therefore moments of $\tilde{\beta}_\pi$ and related quantities exist.

The purpose of the experimenter is to minimize the error $E_\pi \int_{-1}^1 \frac{1}{2}\{ \tilde{\ell}_\pi(x) - \ell(x) \}^2 {dx}$. However, this quantity is difficult to handle because the dependence on $\pi$ can be complicated. Therefore, Theorem \ref{thm:bound} below states an asymptotic approximation where the dependence on $\pi$ is explicit, and leads to the asymptotic minimax design stated in Theorem \ref{thm:main}. To this end, define 
\begin{equation*}
	{\Omega}_\pi= {\rm var}_\pi\left\{\vec{{X}}_1^T e_1 \frac{1/2}{\pi(X_1)} \right\}.
\end{equation*}
The matrix ${\Omega}_\pi$ can be written explicitly using the law of total variance 
\begin{align} \nonumber
	{\rm var}_\pi\left\{\vec{{X}}_1^T e_1 \frac{1/2}{\pi(X_1)} \right\} &=
	E_\pi\left\{ \vec{{X}}_1^T \vec{{X}}_1  \frac{1/4}{\pi^2(X_1)} {\rm var} \left( e_1 \Big | X_1 \right) \right\}
	+{\rm var}_\pi\left\{\vec{{X}}_1^T \frac{1/2}{\pi(X_1)} E ( e_1   | X_1 ) \right\}\\
	\nonumber
	&=\frac{\sigma^2}{4}\int_{-1}^1 \vec{x}^T \vec{x} \frac{1}{\pi(x)} dx+
	\frac{1}{4}\int_{-1}^1 \vec{x}^T \vec{x} \frac{\{m(x)-\ell(x)\}^2}{\pi(x)} dx\\ \label{eq:Omega_def1}
	&=\frac{1}{4} \int_{-1}^1 \vec{x}^T \vec{x} \frac{1}{\pi(x)} [\sigma^2 + \{m(x)-\ell(x)\}^2]dx.
\end{align}
The following theorem quantifies the approximation, $E_\pi \int_{-1}^1 \frac{1}{2}\{ \tilde{\ell}_\pi(x) - \ell(x) \}^2 {dx} \approx \frac{1}{n}tr( {Q}^{-1} {\Omega}_\pi)$.
\begin{theorem}\label{thm:bound}
	Consider model \eqref{eq:model}, and suppose that $E_\pi(e_1^4)$ is finite and $\pi(x)>\ve_0$ for all $x\in[-1,1]$. Then for all $n$,
	\begin{equation}\label{eq:thm_bound}
		n^{3/2} \left|  \int_{-1}^1 \frac{1}{2} E_\pi\{ \tilde{\ell}_\pi(x) - \ell(x) \}^2  dx - \frac{1}{n}tr( {Q}^{-1} {\Omega}_\pi)  \right| < C,
	\end{equation}
	where $C$ is a constant that depends on $\pi$ only through $\ve_0$.
\end{theorem}
Theorem \ref{thm:bound} implies that 
\begin{equation}\label{eq:O_3_2}
	\int_{-1}^1 \frac{1}{2} E_\pi\{ \tilde{\ell}_\pi(x) - \ell(x) \}^2  dx = \frac{1}{n}tr( {Q}^{-1} {\Omega}_\pi) +O(1/n^{3/2}), 
\end{equation}
and the error term is uniform for all designs $\pi$, with $\pi(x)>\ve_0$.  
{It follows that the approximation in \eqref{eq:O_3_2} is uniformly valid for all the densities $\pi$ that are bounded away from zero. Since when $\pi(x)$ is close to zero for some $x$, the minimax risk can be large, it is enough to consider only densities that are bounded away from zero. }

The focus from now on will be on the asymptotic variance term $tr( {Q}^{-1} {\Omega}_\pi)$.
Define the set  
\[
{\cal M}=\left\{ m(x) ~\text{such that}~\int_{-1}^1\frac{1}{2} \{m(x)-\ell(x)\}^2 dx \le 1 \right\},
\]
which bounds the possible deviation from linearity. The bound 1 is arbitrary, and in general the minimax design depends on the ratio between $\sigma^2$ and the linearity bound; see \eqref{eq:R} below and also the discussion in Section \ref{sec:var}.
Let $h(x)= \vec{x}^T {Q}^{-1}  \vec{x}$.
Theorem 3.3 of \citet{Wiens2000} connects $h(x)$ to Legendre polynomials. Let 
\begin{equation}\label{eq:R}
	{\cal R}_\pi=\frac{\sigma^2}{2} \int_{-1}^1 \frac{h(x)}{\pi(x)} dx + \sup_{x \in [-1,1]} \frac{h(x)}{\pi(x)}.
\end{equation}
{Lemma 1 of \citet{Nie2018} implies the following proposition.}
\begin{proposition}\label{prop:sup}
	Under the above setting, $\sup_{m \in {\cal M}} tr( {Q}^{-1} {\Omega}_\pi)= {\cal R}_\pi/2$.
\end{proposition}  
Thus, the minimax problem is to find a density $\pi(x)$ that minimizes ${\cal R}_\pi$.

\section{Main result}\label{sec:main}

Consider the problem of minimizing ${\cal R}_\pi$, which is given in \eqref{eq:R}. The first part is minimized by $\pi(x) \propto h^{1/2}(x)$ and the second part is optimized when $\pi(x) \propto h(x)$; see Remark 3.2 of \citet{Nie2018}. Theorem \ref{thm:main} below shows that the minimizer of ${\cal R}_\pi$ is composed from both solutions. 


To state the result about the minimax design, some definitions are needed. Define the space of designs, which are densities supported in $[-1,1]$,
\[
\Pi = \left \{ \pi: [-1,1] \mapsto \mathbb{R}^+ \text{ such that } \int_{-1}^1 \pi(x)dx=1 \right \}.
\]
The purpose is to find ${\pi}^* \in \Pi$ that minimizes ${\cal R}_{\pi}$ for $\pi \in {\Pi}$. 
Consider $h_0 \in [h_{min},h_{max}]$, where $h_{min} = \min_{x \in [-1,1]} h(x)$ and $h_{max} = \max_{x \in [-1,1]} h(x)$. Define, 
\begin{equation}\label{eq:A_h0}
	A_{h_0}=\{ x \in [-1,1] \text { such that } h(x) \le h_0\} 
	\text{ and }B_{h_0}=\{ x \in [-1,1] \text { such that } h(x) > h_0\}.
\end{equation} 
The sets $A_{h_0}$ and $B_{h_0}$ have an explicit expression.
Let $x_{0,1},\ldots,x_{0,J}$ be all the solutions of $h(x)=h_0$ for $x \in [-1,1]$ ordered from smallest to largest. If none of $x_{0,1},\ldots,x_{0,J}$ is an extremal point of $h$, then $A_{h_0}= [x_{0,1},x_{0,2}] \cup [x_{0,3},x_{0,4}] \ldots  \cup [x_{0,J-1},x_{0,J}]$ and $B_{h_0}= [-1,x_{0,1}] \cup [x_{0,2},x_{0,3}] \ldots  \cup [x_{0,J},1]$.
Define the function
\begin{equation}\label{eq:f_h0}
	f(h_0)=  \frac{\int_{B_{h_0}}\{h_0- h(x)\}dx}{h_0}.
\end{equation}
Since  $h(x) \ge h_0$ for $x \in B_{h_0}$, we have that $f(h_0) \le 0$; also, $f(h_{max})=0$. A calculation shows that $f'(h_0)= \frac{\int_{B_{h_0}} h(x)dx}{h_0^2} >0$, i.e., $f(h_0)$ is monotone increasing. {Let $\sigma^2_{min}=-\frac{2}{f(h_{min})}$. } The main result is given now.  

\begin{theorem}\label{thm:main}
	Under the above setting,
	\begin{enumerate}[(i)]
		\item If $\sigma^2\le \sigma^2_{min}$, the design ${\pi}^*(x) = \frac{h(x)}{\int_{-1}^1 h(x) dx }$ satisfies  ${\cal R}_{{\pi}^*} \le {\cal R}_{\pi}$ for all $\pi \in {\Pi}$.
		\item If $\sigma^2 > \sigma^2_{min}$, define $h_0^*$ to be such that $f(h_0^*)=-\frac{2}{\sigma^2}$.  The design
		\begin{equation}\label{eq:form3}
			{\pi}^*(x)=\left\{ \begin{array}{cc} c\{h^*_0 h(x)\}^{1/2} &  x \in A_{h_0^*} \\  
				c h(x) &  x \in B_{h_0^*} \end{array} \right. ,
		\end{equation}	
		where $c$ is a normalizing constant, satisfies  ${\cal R}_{\pi^*} \le {\cal R}_{\pi}$ for all $\pi \in {\Pi}$.
	\end{enumerate}	
\end{theorem}

By the properties of the function $f$ that were mentioned above, when $\sigma^2 > \sigma^2_{min}$ the solution $f(h_0^*)=-\frac{2}{\sigma^2}$ exists and is unique. Since $h(x_{0,j})=h_0$, the density ${\pi}^*(x)$ is well-defined for $x=x_{0,j}$, $j=1,\ldots,J$, and ${\pi}^*(x)$ is continuous. 

The critical value $\sigma^2_{min}$ is equal to $\frac{2}{\int_{-1}^1 h(x)dx/h_{min} -2}$. For $K=1$, $\sigma^2_{min}=1$ and for $K=2$, $\sigma^2_{min}=3/2$. When $K$ is large, $h(x) \approx \frac{2 K}{ (1-x^2)^{1/2}  \beta^2(3/4,3/4) }$ \citep[][Theorem 3.3]{Wiens2000} and  $\sigma^2_{min}\approx \frac{2}{\pi -2} \approx 1.752$.    

Minimization of ${\cal R}_\pi$ is similar to the problem that is studied in \citet{Wiens1998} and \citet{Wiens2000}. The same minimax criterion is considered for the weighted least squares estimators. However, these two papers study the estimation error when $X_1,\ldots,X_n$ are conditioned upon. Specifically, they consider minimization of the first part of ${\cal R}_\pi$ as defined in \eqref{eq:R}, which is an asymptotic approximation of $\sup_{m \in {\cal M}}\left[ E_\pi \left\{ \int_{-1}^1 {\rm var}( \tilde{\ell}_\pi(x)|X_1,\ldots X_n) dx\right\}\right]$. This part pertains to the variance of $\tilde{\ell}_\pi(x)$ conditioned on the design, i.e., when the design is considered fixed; see \eqref{eq:Omega_def1}. \citet{Wiens2000} shows that it is minimized when $\pi(x) \propto \sqrt{h(x)}$. 
On the other hand, the second part of ${\cal R}_\pi$ comes from the variability of $E( \tilde{\ell}_\pi(x)|X_1,\ldots X_n)$ and is not considered in \citet{Wiens1998} and \citet{Wiens2000}.

As mentioned in the introduction, \citet{Nie2018} study the same problem of minimizing ${\cal R}_\pi$ in the active learning setup. In particular, Theorem 1 of \citet{Nie2018} shows that in the current setting the optimal design has the form $\pi^*(x)=\left\{ \begin{array}{cc} c\{h_0^* h(x)\}^{1/2} &  x \in A_{h^*_0} \\  
	c h(x) &  x \in B_{h^*_0} \end{array} \right. $ for some optimal $h^*_0$. This optimal form applies to both the discrete and continuous settings. Theorem 2 of the current work shows that in the continuous case the optimal $h^*_0$ must satisfy $f(h^*_0)=-2/\sigma^2$ when such a solution exists and otherwise it equals to $h_{min}$. It is also shown that $f$ is monotone, and therefore the solution if exists is unique and can be easily found using numerical algorithms. This result implies the phase-transition phenomenon of Theorem \ref{thm:main}, namely, that when $\sigma^2 \le \sigma^2_{min}$ the optimal design is ${\pi}^*(x) = \frac{h(x)}{\int_{-1}^1 h(x) dx }$ and for $\sigma^2> \sigma^2_{min}$ the optimal design is given in \eqref{eq:form3}.



\section{Unknown variance and Heteroscedasticity}\label{sec:var}

The results so far assumed that the conditional variance $\sigma^2$ is known. When this is not the case, one can consider a parameter space where $\sigma^2$ is upper bounded by $\bar{\sigma}^2$, say. In this case the equivalent result for Proposition \ref{prop:sup}, with essentially the same proof, is
$\sup_{m \in {\cal M}, \sigma^2 \le \bar{\sigma}^2 } tr( {Q}^{-1} {\Omega}_\pi)= {\cal R}_\pi/2$,
where, ${\cal R}_\pi=\frac{\bar{\sigma}^2}{2} \int_{-1}^1 \frac{ h(x)}{\pi(x)} dx + \sup_{x \in [-1,1]} \frac{h(x)}{\pi(x)}$.
Then Theorem \ref{thm:main} applies with $\sigma^2=\bar{\sigma}^2$. Thus, the minimax design depends on the ratio between the upper bound on the variance, $\bar{\sigma}^2$, and the bound on the departure from the linear model quantified by  $\int_{-1}^1 \frac{1}{2} \left\{ m(x)-\ell(x)\right\}^2 dx$, which was arbitrarily set to one. Practically, the experimenter can choose this ratio based on the specific problem at hand and compute the minimax design with respect to the selected ratio. When this ratio is smaller or equal to $\sigma^2_{min}$, the minimax design is the same. In particular, for $K=1$ and $K=2$ if this ratio is one, i.e., the bound of the variance is equal to the bound of the departure from linearity, the minimax design is the same as when the bound of the variance is zero.   

Model \eqref{eq:model} assumes that ${\rm var}(Y_i | X_i)$ is a constant number, $\sigma^2$, which does not depend on $X_i$. This is a strong assumption, and the results for the heteroscedastic case are now discussed. Suppose that ${\rm var}(Y_i|X_i)=\sigma^2(X_i)$; here $\sigma^2$ is a function of $X$, rather than a constant. The  estimator $\tilde{\ell}_\pi(x)$, which is defined after \eqref{eq:calE}, is still a consistent estimator of $\ell(x)$, but now the equivalent expression for ${\Omega}_\pi$ instead of \eqref{eq:Omega_def1} is
${\Omega}_\pi= \int_{-1}^1 \vec{x}^T \vec{x} \frac{1/4}{\pi(x)} [\sigma^2(x) + \{m(x)-\ell(x)\}^2]dx$.
Therefore, ${\cal R}_{\pi}$ changes to
${\cal R}_\pi=\frac{1}{2} \int_{-1}^1 \frac{ \sigma^2(x) h(x)}{\pi(x)} dx + \sup_{x \in [-1,1]} \frac{h(x)}{\pi(x)}$,
i.e., $\sigma^2(x)$ replaces $\sigma^2$. 
The minimax design ${\pi}^*(x)$, thus, depends on $\sigma^2(x)$, and is not given here. However, the function $\sigma^2(x)$ is rarely known in advance and instead it is natural to consider a setting where $\sigma^2(x) \le \bar{\sigma}^2$ and then Theorem \ref{thm:main} with $\sigma^2=\bar{\sigma}^2$ applies as above. The situation where $\sigma^2(x)$ is bounded by a function $\bar{\sigma}^2(x)$ is deferred to future work.




\section{Numerical illustrations}\label{sec:numerical}

In this section the optimal minimax design, ${\pi}^*$, which is given in Theorem 2, is illustrated in a few examples and a comparison to other designs is presented.

\begin{figure}[ht!]
	\subfigure[$K=1$]{%
		\includegraphics[width=0.45\textwidth]{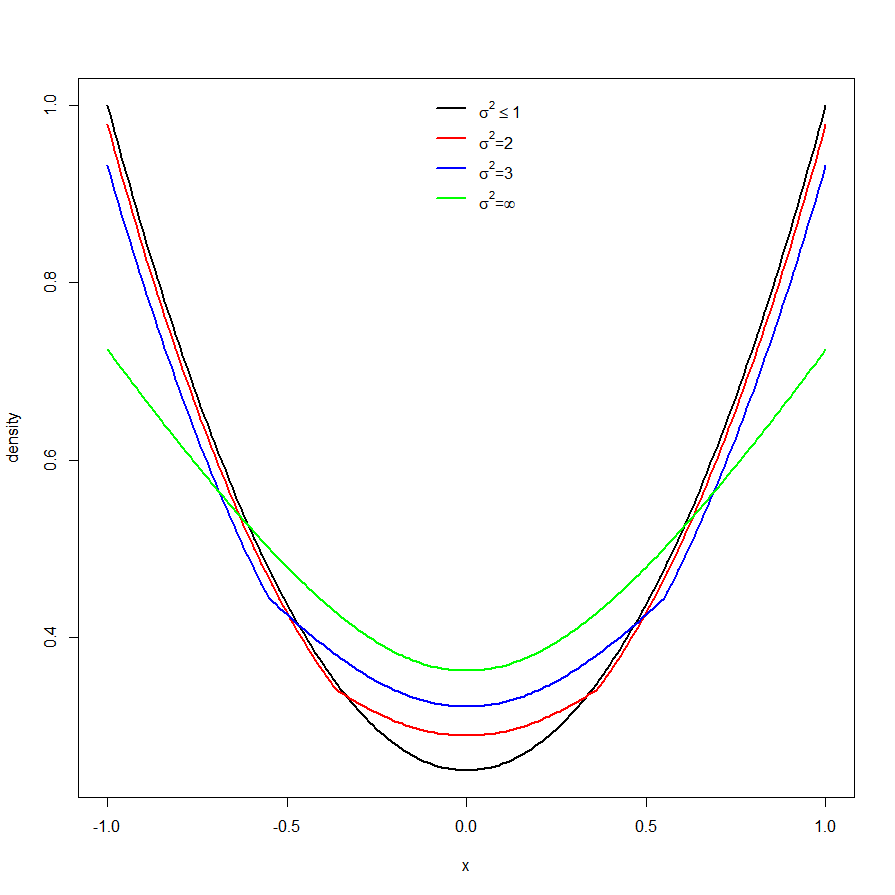}
	}
	\subfigure[$K=2$]{%
		\includegraphics[width=0.45\textwidth]{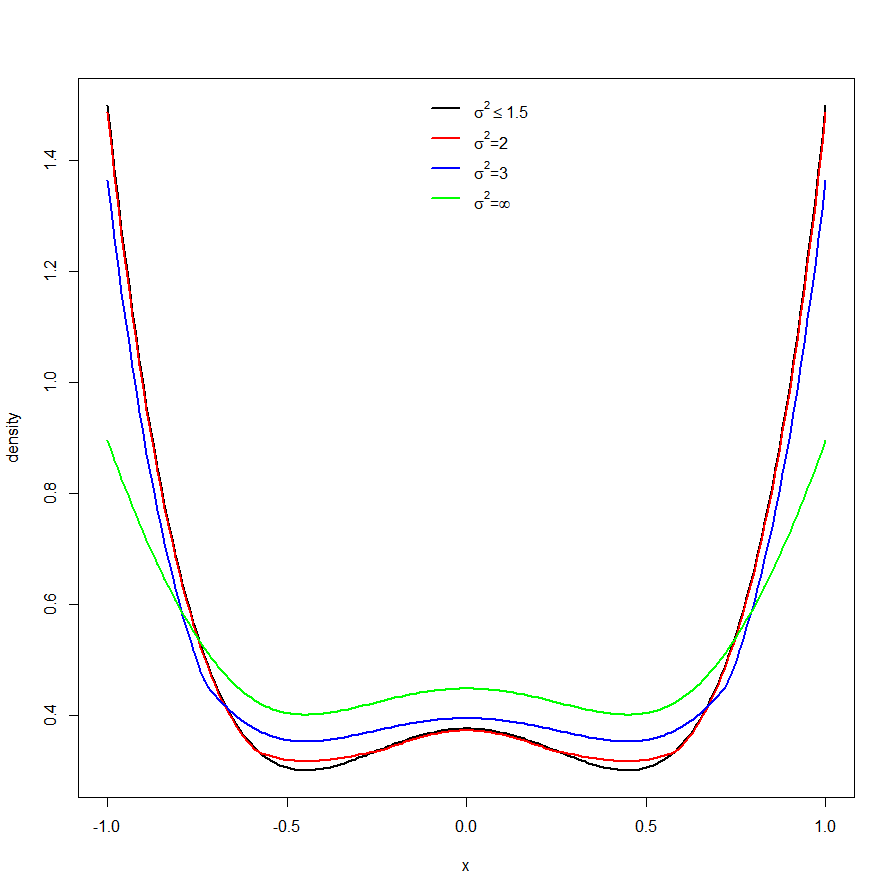}
	}	
	\caption{The optimal design ${\pi}^*$ for $\sigma^2 \le \sigma^2_{min}$ , $\sigma^2=2,3,\infty$ and $K=1,2$.}\label{fig:K_1_2}
\end{figure}

Figure \ref{fig:K_1_2} plots the optimal design ${\pi}^*$ for $\sigma^2 \le \sigma^2_{min}$, $\sigma^2=2,3,\infty$ and $K=1,2$. When $\sigma^2 \le \sigma^2_{min}$, the optimal design is proportional to $h(x)$. When $\sigma^2=2,3$, for some $x$'s, ${\pi}^*(x)$ is proportional to $h^{1/2}(x)$ and for others, it is proportional to $h(x)$. The intervals where the former case holds are given in Table \ref{tab:1}. For $\sigma^2=\infty$, ${\pi}^*(x)$ is proportional to $h^{1/2}(x)$ for all $x$.

\begin{table} [ht!]
	\begin{center}
		\begin{tabular}{c||c|c}
			{} & $K=1$ & $K=2$\\
			\hline 
			\hline
			$\sigma^2 \le \sigma^2_{min}$  &  $\emptyset$ &  $\emptyset$ \\
			$\sigma^2=2$  &  $[-0.364,0.364]$ & $[-0.587,-0.235] \cup [0.235,0.587]$\\
			$\sigma^2=3$ &  $[-0.550,  0.550]$ & $[-0.725,  0.725]$\\
			$\sigma^2=\infty$ & $[-1,1]$ & $[-1,1]$
	\end{tabular}
\end{center}
\caption{The intervals $A_{h_0^*}$}
\label{tab:1}
\end{table}


As Huber writes \citep[][page 244]{huber2004robust} ``with any minimax
procedure, there is the question of whether it is too pessimistic and perhaps safeguards
only against some very unlikely contingency''. In order to examine this question, the minimax design is compared against the sqrt design of \citet{Wiens2000}, which corresponds to $\sigma^2=\infty$ in the above setting, and the uniform design, in scenarios where $m(x)$ is a quadratic function of $x$ for $K=1$ or a cubic function for $K=2$. In the former case, $m(x)=x+3.354 x^2$ and in the latter case $m(x)=x+x^2/2 + 6.614 x^3$ were examined. The leading coefficient of $m(x)$ in both cases were set to satisfy the condition $\int_{-1}^1\frac{1}{2} \{m(x)-\ell(x)\}^2 dx = 1$. It is also of interest to investigate the small sample performance of the minimax design. 

Table \ref{tab:2} presents empirical evaluations of $n \times \int_{-1}^1 \frac{1}{2} E_\pi\{ \tilde{\ell}_\pi(x) -\ell(x) \}^2 dx$ for $n=50$, which are based on $10^5$ simulations. Model (3)  was used  with normal $\ve$'s to generate the simulated dataset. Also presented is the asymptotic approximation of the above expression, which is
\[
tr({Q}^{-1} {\Omega}_\pi)=\frac{\sigma^2}{4} \int_{-1}^{1} \frac{h(x)}{\pi(x)} dx+\frac{1}{4}\int_{-1}^{1} \frac{h(x)}{\pi(x)} \{m(x)- \ell(x)\}^2 dx.
\] 
In order to better compare the designs, for each simulated dataset, $U_1,\ldots,U_{50} \sim U(0,1)$ were sampled and $X_i$ was defined to be $2 U_i -1, F_{sqrt}^{-1} (U_i), F_{minimax}^{-1} (U_i)$ for the uniform, sqrt, minimax designs respectively, where $F_{sqrt}$ and $F_{minimax}$ are the cumulative distribution functions of the sqrt and minimax densities, $i=1,\ldots,50$. Furthermore, the same $\ve$'s were used for all three designs. It follows that the standard error of the differences of the empirical means are smaller than the standard errors of each mean.  The former is given in the second last column of Table \ref{tab:2} for the differences of the means of the minimax and sqrt designs,  and the latter is given after each mean in parentheses. The best design for each setting is marked in italic font, and when the difference is not statistically significant, i.e., less than two standard errors, also the second best is marked. The last column of Table \ref{tab:2} shows the percentage of improvement between the best and the second best designs.         

In all scenarios, the uniform design has the largest estimation error and the ensuing discussion focuses on the comparison between the sqrt and minimax designs. 
When $\sigma^2=1,1/2$ the best design is the minimax and the improvement with respect to the sqrt design is higher when $n=50$ compared to the asymptotic scenario, when $K=2$ compared to $K=1$ and when $\sigma^2=1/2$ compared to $\sigma^2=1$. For $\sigma^2=2,3$ the sqrt design is generally better than minimax, but the percentage of improvement is smaller compared to $\sigma^2=1,1/2$ and in most cases the estimation error is almost the same. The conclusion is that it seems safe to use the minimax designs under reasonable models for $m(x)$ and for relatively small sample sizes.   

\begin{table} [ht!]
{\small
		\begin{tabular}{c|c|c||c|c|c||c|c| }
			{}&{}&{}&uniform&sqrt&minimax& SE & \% improvement\\
			\hline \hline
			\multirow{8}{*}{$K=1$} & \multirow{2}{*}{$\sigma^2=1/2$} & $n=50$ & $3.42 (0.04)$ & $2.85 (0.03)$ & ${\it 2.69 (0.03)}$ & $0.01$ & $5.6\%$\\
			& & asymptotic & $3.07$ & $2.62$ & ${\it 2.50}$ &{} & $4.6\% $\\ 
			\cline{2-8}
			& \multirow{2}{*}{$\sigma^2=1$} & $n=50$ & $4.93 (0.05)$ & $4.31 (0.04)$ & ${\it 4.22 (0.04)}$ & $0.02$ & $2.1 \%$ \\
			& & asymptotic & $4.57$ & $4.04$ & ${\it 4.00}$ &  &$0.99\%$ \\
			\cline{2-8}
			& \multirow{2}{*}{$\sigma^2=2$} & $n=50$ & $11.21 (0.12)$ & ${\it 10.29 (0.11)}$ & ${\it 10.31 (0.11)} $& $0.03$ & $0.2 \% $\\
			& & asymptotic & $10.57$ & ${\it 9.76}$ & $9.84$ & & $0.8 \%$ \\
			\cline{2-8}
			& \multirow{2}{*}{$\sigma^2=3$} & $n=50$ & $20.95 (0.21)$ & ${\it 19.55 (0.20)}$ & ${\it 19.62 (0.20)}$ & $0.04$ & $0.36\%$ \\
			& & asymptotic & $20.57$ & ${\it 19.29} $& $19.38$ & & $0.46 \%$ \\
			\hline
			\multirow{8}{*}{$K=2$} & \multirow{2}{*}{$\sigma^2=1/2$} & $n=50$ & $6.21 (0.08)$ & $4.49 (0.05)$ & ${\it 4.06 (0.04)}$ & $0.03$ & $9.58 \%$ \\
			& & asymptotic & $5.13$ & $4.03$ & ${\it 3.72}$ & & $7.7 \%$  \\ 
			\cline{2-8}
			& \multirow{2}{*}{$\sigma^2=1$} & $n=50$ & $8.77 (0.10)$ & $6.76 (0.06)$ & ${\it 6.42 (0.06)}$ & $0.04$ & $5.03\%$ \\
			& & asymptotic & $7.38$ & $6.13$ & ${\it 5.98}$ & &$2.45 \%$ \\
			\cline{2-8}
			& \multirow{2}{*}{$\sigma^2=2$} & $n=50$ & $18.46 (0.18)$ & ${\it 15.61 (0.14)}$ & ${\it 15.67 (0.14)}$ & $0.07$ & $0.38 \%$ \\
			& & asymptotic & $16.38$ & ${\it 14.57}$ & $14.92$ & &$3.02\%$ \\
			\cline{2-8}
			& \multirow{2}{*}{$\sigma^2=3$} & $n=50$ & $34.37 (0.32)$ & ${\it 30.34 (0.27)}$ & ${\it 30.48 (0.27)}$ & $0.08$ & $0.46 \%$ \\
			& & asymptotic & $29.84$ & ${\it 28.62}$ & $29.02$ & &$1.38\% $ \\
			\hline
	\end{tabular}}
	\label{tab:2}	
		\caption{\footnotesize Simulation means  (standard errors) of  $n \times \int_{-1}^1 \frac{1}{2} E_\pi\{ \tilde{\ell}_\pi(x) -\ell(x) \}^2 dx$ when $n=50$ and the asymptotic approximation $tr({Q}^{-1} {\Omega}_\pi)$ for the uniform, sqrt and minimax designs, $\sigma^2=1/2,1,2,3$ and $K=1,2$ are presented. For each setting, the best design, and the second best if the difference is not statistically significant, is marked in italic font. The second last column shows the simulation standard error of the difference between the means of the minimax and sqrt designs, and the last column shows the percentage of improvement between the best and the second best designs.}
\end{table}

\section{General covariates} \label{sec:gen}

So far the discussion was limited to the case where the support of $X_i$ is $[-1,1]$, where the estimation error of $\tilde{\ell}_\pi(x)$ was equally weighted in $x$ and where the covariates were monomials of $x$. The general setting is presented now. 

Suppose that the support of $X_i$ is a compact set ${S} \subset {\mathbb R}^p$ and the estimation error is weighted by $\lambda(x)$, which is a density on ${S}$ with respect to the Lebesgue measure. It is assumed that $\lambda(x)$ is a continuous function of $x$ and $\lambda(x) >0$ for all $x \in {S}$. The conditional expectation $m(x)=E(Y_i|X_i=x)$ is approximated by a vector of covariates $\vec{x}=(c_1(x),c_2(x),\ldots,c_{k}(x))$ for some predetermined continuous functions $c_1,\ldots,c_K$, and the best linear approximation is $\ell(x)=\vec{x} {\beta}$, where ${\beta}={Q}^{-1} \int_{S} \vec{x}^T m(x) \lambda(x) dx$ and ${Q}=\int_{S} \vec{x}^T  \vec{x} \lambda(x) dx$.  
The experimenter selects a density $\pi$ on ${S}$ and then samples $(X_1,Y_1),\ldots,(X_n,Y_n)$, where $X_i \sim \pi$. The goal is to find $\pi \in \Pi$, where
\[
\Pi = \left \{ \pi: {S} \mapsto \mathbb{R}^+ \text{ such that } \int_{S} \pi(x)dx=1 \right \},
\]
that is optimal in a minimax sense. 

In order to estimate ${\beta}$, the same weighted least squares estimator is used as in Section 2, but now the weights are $\frac{\lambda(X_i)}{\pi(X_i)}$. This leads to the estimator $\tilde{\ell}_\pi(x)$ as in Section 2. By the same arguments of Theorem 1, we have
\[
\int_{S} E_\pi\{ \tilde{\ell}_\pi(x) - \ell(x) \}^2 \lambda(x) dx = \frac{1}{n}tr( {Q}^{-1} {\Omega}_\pi) +O(1/n^{3/2}), 
\]
where,
\[
{\Omega}_\pi= \int_{S} \vec{x}^T \vec{x} \frac{\lambda^2(x)}{\pi(x)} [\sigma^2 + \{m(x)-\ell(x)\}^2]dx.
\] 

Consider the set
\[
{\cal M}=\left\{ m(x) ~\text{such that}~\int_{S}  \{m(x)-\ell(x)\}^2  dx \le 2 \right\};
\]
for technical reasons it is convenient to weigh $\{m(x)-\ell(x)\}^2$ equally, rather than by $\lambda(x)$ and to set the bound to two.
Then, similar to Proposition 1, we have 
$\sup_{m \in {\cal M}} tr( {Q}^{-1} {\Omega}_\pi)=  {\cal R}_\pi/2$, where,
\[
{\cal R}_\pi=\frac{\sigma^2}{2} \int_{S} \frac{h(x) }{\pi(x)}  dx + \sup_{x \in [-1,1]} \frac{h(x)}{\pi(x)},
\]
and $h(x)= 4 \vec{x}^T {Q}^{-1}  \vec{x} \lambda^2(x)$.

As in Section 3, let $h_{min}=\min_{x \in {S}} h(x)$ and $h_{max}=\max_{x \in {S}} h(x)$. Consider $h_0 \in [h_{min},h_{max}]$ and define the sets
$A_{h_0}= \{x \in {S} ~:~h(x) \le {h_0}  \}$, $B_{h_0}= \{x \in {S} ~:~h(x) > {h_0}\}$ and the function $f(h_0)=  \frac{\int_{B_{h_0}}\{h_0- h(x)\}dx}{h_0}$. {Let also $\sigma^2_{min}=-\frac{2}{f(h_{min})}$. } Then Theorem 2 applies with this notation using similar arguments, i.e., if $\sigma^2\le \sigma^2_{min}$, the design ${\pi}^*(x)= \frac{h(x)}{\int_{{S}} h(x) dx }$ satisfies  ${\cal R}_{{\pi}^*} \le {\cal R}_{\pi}$ for all $\pi \in {\Pi}$, and if $\sigma^2> \sigma^2_{min}$ the minimax design ${\pi}^*$ is given in (12).

\section{Proofs}
\subsection*{Proof of Theorem 1}

Write
\begin{multline*}
	E_\pi \left[ \int_{-1}^1  \frac{1}{2} \{ \tilde{\ell}_\pi(x) - \ell(x) \}^2 dx\right]= tr \left[ E_\pi\left\{(\tilde{\beta}_\pi - {\beta}) (\tilde{\beta}_\pi - {\beta})^T  \right\} \int_{-1}^1  \frac{1}{2}\vec{x}^T \vec{x}dx\right]\\
	= tr \left[ E_\pi\left\{(\tilde{\beta}_\pi - {\beta}) (\tilde{\beta}_\pi -{\beta})^T  \right\} {Q} \right].
\end{multline*}
By Theorem 5.1 of \citet{Tropp}, ${\rm pr}({\cal E}_\pi)\le a^n$ where $a<1$, because $E_\pi[\tilde{Q}_\pi]=n{Q}$, and since when ${\cal E}_\pi$ occurs $\tilde{\beta}_\pi$ is bounded, we have that
\begin{equation}\label{eq:tr_G_xi}
	tr \left[ E_\pi\left\{(\tilde{\beta}_\pi -{\beta}) (\tilde{\beta}_\pi - {\beta})^T  \right\} {Q} \right]=tr \left[ E_\pi\left\{(\tilde{\beta}_\pi - {\beta}) (\tilde{\beta}_\pi - {\beta})^T   I ({\cal E}_\pi^C)\right\} {Q} \right] + b_n,
\end{equation}
where $b_n$ satisfies $|b_n| \le b^n$ for $b<1$. When ${\cal E}_\pi^C$ occurs, then
\begin{multline*}
	\tilde{\beta}_\pi -{\beta} =\left( {\mathbb{X}}^T {W}_\pi {\mathbb{X}} \right)^{-1}{\mathbb{X}}^T {W}_\pi {\mathbb{Y}} - {\beta}= \left( {\mathbb{X}}^T {W}_\pi {\mathbb{X}} \right)^{-1}{\mathbb{X}}^T {W}_\pi ({\mathbb{X}} {\beta} + {\mathbb E} ) - {\beta}\\
	=\left( {\mathbb{X}}^T {W}_\pi {\mathbb{X}} \right)^{-1}{\mathbb{X}}^T {W}_\pi {\mathbb E}=
	\frac{1}{\sqrt{n}}\left(\frac{1}{n} \tilde{Q}_\pi \right)^{-1} \frac{1}{\sqrt n}{\mathbb{X}}^T {W}_\pi {\mathbb E} = \frac{1}{\sqrt{n}} \bar{Q}_\pi^{-1} \bar{U}_\pi, 
\end{multline*} 
where ${\mathbb E}=(e_1,\ldots,e_n)^T$, $\bar{Q}_\pi=\tilde{Q}_\pi/n$, $\bar{U}_\pi=\frac{1}{\sqrt{n}}\sum_{i=1}^n {U}_{i,\pi}$ and ${U}_{i,\pi}= \vec{{X}}_i^T \frac{1/2}{\pi(X_i)} e_i$.  We have that $E_\pi({U}_{i,\pi})={ 0}$, $i=1,\ldots,n$, because
\begin{multline*}
	E_\pi\left( {U}_{1,\pi} \right)= E_\pi\left\{ \vec{{X}}_1^T \frac{1/2}{\pi(X_1)} e_1 \right\}=  \int_{-1}^1 \frac{1/2}{\pi(x)}\vec{x}^T E(e_1|X_1=x) \pi(x) dx\\
	=\int_{-1}^1 \frac{1}{2}\vec{x}^T E(e_1|X_1=x)  dx=\int_{-1}^1 \frac{1}{2}\vec{x}^T \{ m(x)-\ell(x)\}  dx ={ 0},
\end{multline*}
where the last equality follows from the orthogonality property of $\ell(x)$.
Also, it is easy to verify that $E_\pi({U}_{1,\pi} {U}_{1,\pi}^T) = {\Omega}_\pi$. 

Going back to \eqref{eq:tr_G_xi}, 
\[
tr \left[ E_\pi\left\{(\tilde{\beta}_\pi - {\beta}) (\tilde{\beta}_\pi -{\beta})^T   I ({\cal E}_\pi^C)\right\} {Q} \right] = \frac{1}{n}
tr \left[ E_\pi\left\{ \bar{Q}_\pi^{-1} \bar{U}_\pi \bar{U}_\pi^T  \bar{Q}_\pi^{-1} I ({\cal E}_\pi^C)\right\} {Q} \right].
\]
Therefore,
\begin{multline*}
	n \left\{ E_\pi\left[ \int_{-1}^1 \frac{1}{2}\{ \tilde{\ell}_\pi(x) - \ell(x) \}^2dx\right] - tr( {Q}^{-1} {\Omega}_\pi)\right\}\\
	= tr \left( \left[ E_\pi\left\{ \bar{Q}_\pi^{-1} \bar{U}_\pi \bar{U}_\pi^T  \bar{Q}_\pi^{-1} I ({\cal E}_\pi^C)\right\}  -{Q}^{-1} {\Omega}_\pi {Q}^{-1}  \right] {Q}\right)+c_n,
\end{multline*}
where $c_n$ satisfies $|c_n|\le c^n$ for $c<1$.
Thus, in order to show (7), the latter trace needs to be bounded. To this end, write
\begin{multline}\label{eq:I_II_III}
	E_\pi\left\{ \bar{Q}_\pi^{-1} \bar{U}_\pi \bar{U}_\pi^T  \bar{Q}_\pi^{-1} I ({\cal E}_\pi^C)\right\}  -{Q}^{-1} {\Omega}_\pi {Q}^{-1}  =  E_\pi\left\{ \left(\bar{Q}_\pi^{-1} -  {Q}^{-1}\right) \bar{U}_\pi \bar{U}_\pi^T  \bar{Q}_\pi^{-1} I ({\cal E}_\pi^C)\right\} \\
	+ {Q}^{-1} E_\pi\left\{  \bar{U}_\pi \bar{U}_\pi^T  \left(\bar{Q}_\pi^{-1}- {Q}^{-1} \right) I ({\cal E}_\pi^C)\right\} + {Q}^{-1} E_\pi \left\{ \bar{U}_\pi \bar{U}_\pi^T  I ({\cal E}_\pi^C) -  {\Omega}_\pi \right\}  {Q}^{-1}\\
	=I + II + III. 
\end{multline}
In order to prove the theorem, we need to bound $I,II,III$. We start with bounding $I$. We have
\[
\bar{Q}_\pi^{-1} -  {Q}^{-1}= {Q}^{-1} ( {Q} - \bar{Q}_\pi ) \bar{Q}_\pi^{-1}.
\]
Therefore,
\[
{I}= {Q}^{-1} E_\pi\left\{ \left(\bar{Q}_\pi -  {Q} \right) {Q}^{-1}\bar{U}_\pi \bar{U}_\pi^T  \bar{Q}_\pi^{-1} I ({\cal E}_\pi^C)\right\}.
\]
An element of the matrix $ \left(\bar{Q}_\pi -  {Q} \right) {Q}^{-1}\bar{U}_\pi \bar{U}_\pi^T  \bar{Q}_\pi^{-1} I ({\cal E}_\pi^C)$ is of the form $\sum_j {\cal A}_j {\cal B}_j$, where ${\cal A}_j$ and ${\cal B}_j$ are elements of the matrices $\bar{Q}_\pi -  {Q}$ and ${Q}^{-1}\bar{U}_\pi \bar{U}_\pi^T  \bar{Q}_\pi^{-1} I ({\cal E}_\pi^C)$, respectively. In order to bound $I$, it is enough to consider $E_\pi(| {\cal A}_j {\cal B}_j|)$.
By the Cauchy-Schwarz inequality, $E_\pi(| {\cal A}_j {\cal B}_j|) < \left\{E_\pi \left({\cal A}_j^2\right)\right\}^{1/2} \left\{E_\pi \left({\cal B}_j^2\right)\right\}^{1/2}$. Consider first $E_\pi\left({\cal B}_j^2\right)$; $\bar{Q}_\pi^{-1}I ({\cal E}_\pi^C)$ is bounded because every entry in a positive definite matrix is smaller than the maximal eigenvalue and the latter is bounded by the definition of the event ${\cal E}_\pi$ in (6). Elements of $\bar{U}_\pi \bar{U}_\pi^T$ have second moments because $X_i$ are bounded (and so is $1/\pi(X_i)$) and $e_i$ has a finite second moment. It follows that $E_\pi\left({\cal B}_j^2\right)$ is bounded.  Consider now  $E_\pi \left({\cal A}_j^2\right)$. We have for certain $j_1,j_2$,
\begin{multline*}
	E_\pi \left({\cal A}_j^2\right)= E_\pi\left\{\frac{1}{n} \sum_{i=1}^n \left( \frac{1}{\pi(X_i)} X_i^{j_1} X_i ^{j_2}  - \{{Q} \}_{j_1,j_2}\right)\right\}^2\\
	=\frac{1}{n^2} \sum_{i=1}^n E_\pi \left( \frac{1}{\pi(X_i)} X_i^{j_1} X_i ^{j_2}  - \{{Q} \}_{j_1,j_2}\right)^2 \le C/n,
\end{multline*}  
where the second equality is true because $E_\pi \left( \frac{1}{\pi(X_i)} X_i^{j_1} X_i ^{j_2}\right) = \{{Q} \}_{j_1,j_2}$ and the mixed terms fall. It follows that $|I| \le C/\sqrt{n}$.

Term II in \eqref{eq:I_II_III} can be bounded similarly to I. Consider now Term III. Again, by Theorem 5.1 of \citet{Tropp}, when $e$ have a bounded fourth moment, $ E_\pi \left\{ \bar{U}_\pi \bar{U}_\pi^T  I ({\cal E}_\pi^C)\right\}$ is equal to $E_\pi \left( \bar{U}_\pi \bar{U}_\pi^T \right)$ up to an error term that is bounded by $d^n$ for $d<1$, which is negligible. Therefore, since $ E_\pi \left( \bar{U}_\pi \bar{U}_\pi^T  \right)= {\Omega}_\pi $, the proof of the theorem is now completed. \qed

\subsection*{Proof of Theorem 2} 
For a given $\pi_0 \in {\Pi}$, define $y_0=\sup_{x \in [-1,1]}\frac{h(x)}{\pi_0(x)}$. By Lemma \ref{lem:calc_var} below, every solution of the variation problem
\begin{equation*}
	\text{minimize }\int_{-1}^1 \frac{h(x)}{\pi(x)} dx \text{ under the constraints } \int_{-1}^1 \pi(x) dx = 1 \text{ and }\pi(x) \ge h(x)/y_0,
\end{equation*}
denoted by ${\pi}^*_0$, satisfies that for all $x\in [-1,1]$, either  ${\pi}^*_0(x)=c_0 h^{1/2}(x)$, for a normalizing constant $c_0$, or ${\pi}^*_0(x) = h(x)/y_0$,  and ${\pi}^*_0$ is continuous. It follows that
${\pi}^*_0$ has  the form 
\begin{equation}\label{eq:form1}
	{\pi}^*_0(x)=\left\{ \begin{array}{cc} ch_0^{1/2} h^{1/2}(x) &  x \in A_{h_0} \\  
		c h(x) &  x \in B_{h_0} \end{array} \right. ,
\end{equation}	
where $h_0 = y_0^2 c_0^2$, $c=c_0/h_0^{1/2}$ and $A_{h_0},B_{h_0}$ are defined in (10). Notice that these definitions imply that ${\pi}^*_0$ is continuous. Therefore, for every  $\pi_0 \in {\Pi}$ there exists ${\pi}^*_0$ of the form \eqref{eq:form1} such that ${\cal R}_{{\pi}^*_0} \le {\cal R}_{{\pi}_0}$. {An alternative proof to this result can be found in Theorem 1 of \citet{Nie2018}.}

The problem of minimizing ${\cal R}_{\pi}$ is, thus, reduced to finding the optimal $h_0$ in \eqref{eq:form1}. It can be formalized by
\begin{multline*}
	\min_{h_0}  \frac{1}{c} \left\{ \frac{\frac{\sigma^2}{2}\int_{{A}_{h_0}} h^{1/2}(x)dx}{h_0^{1/2}}  + \frac{\sigma^2}{2} \int_{B_{h_0}} dx  +1 \right\} \\
	\text{ under the constraint } c \left\{\int_{{A}_{h_0}}h^{1/2}(x)dx+ \int_{{B}_{h_0}}{h(x)}dx  \right\}=1.
\end{multline*}
Since under the constraint, $\frac{1}{c}=\int_{{A}_{h_0}}h^{1/2}(x)dx+ \int_{{B}_{h_0}}{h(x)}dx$, one can write the above problem as
\[
\min_{h_0} \left\{\int_{{A}_{h_0}}h^{1/2}(x)dx+ \int_{{B}_{h_0}}{h(x)}dx  \right\}  \left\{ \frac{\int_{{A}_{h_0}} h^{1/2}(x)dx}{h_0^{1/2}}  + \int_{B_{h_0}} dx  +\frac{2}{\sigma^2} \right\}=g(h_0).
\]
Let 
\[
u(h_0)= h_0^{1/2} \int_{A_{h_0}} h^{1/2}(x)dx + \int_{B_{h_0}} h(x)dx ~\text{ and }~v(h_0)=\frac{\int_{A_{h_0}} h^{1/2}(x)dx}{h_0^{1/2}} + {\int_{B_{h_0}} dx+\frac{2}{\sigma^2}};
\]
thus, $g(h_0)=u(h_0) v(h_0)$. We have
\[
u'(h_0)= \frac{1}{2h_0^{1/2}} \int_{A_{h_0}} h^{1/2}(x)dx~\text{ and }~v'(h_0)=\frac{ -\frac{1}{2h_0^{1/2}} \int_{A_{h_0}} h^{1/2}(x)dx}{{h_0}}.
\]
The derivative of $g(h_0)$ is therefore
\begin{multline*}
	g'(h_0)=\\ \frac{1}{2h_0^{1/2}} \int_{A_{h_0}} h^{1/2}(x)dx\left\{\frac{\int_{A_{h_0}} h^{1/2}(x)dx}{h_0^{1/2}} + {\int_{B_{h_0}} dx+ \frac{2}{\sigma^2}} - \frac{\int_{A_{h_0}} h^{1/2}(x)dx}{h_0} - \frac{\int_{B_{h_0}} h(x)dx}{h_0} \right\}\\
	=\frac{1}{2h_0^{1/2}} \int_{A_{h_0}} h^{1/2}(x)dx\left\{ {\int_{B_{h_0}} dx+ \frac{2}{\sigma^2}} - \frac{\int_{B_{h_0}} h(x)dx}{h_0} \right\}\\
	=\frac{1}{2h_0^{1/2}} \int_{A_{h_0}} h^{1/2}(x)dx\left\{ f(h_0)+ \frac{2}{\sigma^2} \right\},
\end{multline*}  
where $f(h_0)$ is defined in (11). Recall the properties of $f$ that were mentioned before Theorem 2. If $\sigma^2\le \sigma^2_{min}$, then $g'(h_0)>0$ for all $h_0$ because $f$ is monotone increasing, and therefore, the optimal $h_0$ is $h_0^*=h_{min}$. If $\sigma^2 > \sigma^2_{min}$, let $h_0^*$ be such that $f(h_0^*)=-\frac{2}{\sigma^2}$. It follows that in this case, $g'(h_0)<0$ (respectively, $g'(h_0)>0$) for $h_0 <  h_0^*$ (respectively, $h_0 >  h_0^*$) and therefore $g(h_0)$ is minimized when $h_0=h_0^*$ and the theorem follows. \qed

\begin{lemma}\label{lem:calc_var}
	Let $\bar{\pi}$ be a minimizer of the variation problem
	\begin{multline}\label{eq:var_problem}
		\text{ minimize } \int_{-1}^1 \frac{h(x)}{\pi(x)} dx \text{ under the constraints}\\
		\int_{-1}^1 \pi(x) dx = 1 \text{ and }\pi(x) \ge k h(x) \text{ for all }x\in[-1,1],
	\end{multline}
	where $k$ is a positive constant. Then for all $x\in [-1,1]$, either  $\bar{\pi}(x)=c h^{1/2}(x)$, for a normalizing constant $c$, or $\bar{\pi}(x) = k h(x)$; furthermore, $\bar{\pi}$ is continuous.  
\end{lemma}	
\noindent{\it Proof of Lemma \ref{lem:calc_var}}\\
Since $\pi(x) \ge k h(x)$, one can write $\pi(x)=\frac{k h(x)}{1+s^2(x)}$ for some function $s(x)$. The variation problem \eqref{eq:var_problem} can be now written in terms of $s(x)$,
\[
\text{ minimize } \int_{-1}^1 s^2(x)dx \text{ under the constraint} \int_{-1}^1 \frac{k h(x) dx}{1+s^2(x)}=1.
\]
By the Lagrange multiplier rule for calculus of variations, the optimal $s^2(x)$ must satisfy
$\{1+s^2(x)\}^2 = \lambda k h(x)$ for some $\lambda$ or equivalently that $s^2(x)=\max\{(\lambda k h(x)-1)^{1/2},0\}$, which implies the result because $\pi(x)=\frac{k h(x)}{1+s^2(x)}$. \qed

\section*{Acknowledgment} I would like to thank Amos Ori for useful discussions about the variation problem of Theorem \ref{thm:main}. The close reading and useful comments of two anonymous referees are also gratefully acknowledged.
 
 \bibliographystyle{chicago}
\bibliography{refs}

\end{document}